\documentclass[runningheads]{llncs}
\usepackage{latexsym}
\usepackage{amssymb}
\usepackage{geompsfi}

\def\Z{{\bf Z}}
\def\C{{\bf C}}
\def\Q{{\bf Q}}
\def\R{{\bf R}}

\def\X{{\rm X}}

\def\SL{{\rm SL}}

\def\ran{{r_{\rm an}}}
\def\MW{Mordell-\mbox{\kern-.16em}Weil}
\font\cyr=wncyr10
\newcommand{\Sha}{\mbox{\cyr X}}

\newcommand{\makefig}[3]{  
  \begin{figure}[htbp]
  \refstepcounter{figure}
  \label{#2}
  \begin{center}
    ~#3~\\
    \medskip
    {\rm Figure \thefigure.  #1}
  \end{center}
  \end{figure}
}

\begin{document}

\title{Curves $Dy^2=x^3-x$ of odd analytic rank}
\titlerunning{Curves $Dy^2=x^3-x$ of odd analytic rank}
\author{Noam D.~Elkies\inst{1}}
\institute{Department of Mathematics, Harvard University,
Cambridge, MA 02138 USA\\{\tt elkies@math.harvard.edu}}
\maketitle

\begin{abstract}
For nonzero rational $D$, which may be taken to be a squarefree integer,
let $E_D$ be the elliptic curve $Dy^2 = x^3 - x$ over~$\Q$
arising in the ``congruent number'' problem.\footnote{
  The problem is: for which~$D$\/ does $E_D$ have nontrivial
  rational points, or equivalently positive rank?
  Such~$D$\/ are called ``congruent''$\!$, because they are precisely
  the numbers that arise as the common difference (``congruum'')
  of a three-term arithmetic progression of rational squares,
  namely the squares of $(x^2-2x-1)/2y$, $(x^2+1)/2y$, and
  $(x^2+2x-1)/2y$.  See the Preface and Chapter~XVI of~\cite{Dickson}
  for the early history of this problem, and~\cite{Koblitz} for a
  more modern treatment of the curves $E_D$.
  }
It is known that the \hbox{$L$-function} of $E_D$
has sign~$-1$, and thus odd analytic rank $\ran(E_D)$,
if and only if $|D|$ is congruent to $5$, $6$, or~$7$ mod~$8$.
For such~$D$, we expect by the conjecture of Birch and Swinnerton-Dyer
that the arithmetic rank of each of these curves $E_D$ is odd,
and therefore positive.
We prove that $E_D$ has positive rank for each $D$\/
such that $|D|$ is in one of the above congruence classes mod~$8$
and also satisfies $|D|<10^6$.  Our proof is computational:
we use the modular parametrization of $E_1$ or $E_2$ to construct
a rational point $P_D$ on each $E_D$ from CM points on modular curves,
and compute $P_D$ to enough accuracy to usually distinguish it from
any of the rational torsion points on $E_D$.
In the $1375$ cases in which we cannot numerically distinguish
$P_D$ from $(E_D)_{\rm tors}$,
we surmise that $P_D$ is in fact a torsion point
but that $E_D$ has rank~$3$, and prove that the rank is positive
by searching for and finding a non-torsion rational point.
We also report on the conjectural extension to $|D|<10^7$
of the list of curves $E_D$ with odd $\ran(E_D)>1$,
which raises several new questions.
\end{abstract}
\section{Introduction}
\subsection{Review: The curves $E_D$ and their arithmetic}
For nonzero rational $D$\/ let $E_D$ be the elliptic curve
\begin{equation}
E_D:  Dy^2 = x^3 - x
\label{E_D}
\end{equation}
over~$\Q$.  Since $E_D$ and $E_{c^2 D}$ are isomorphic
for any nonzero rational $c,D$, we may assume
without loss of generality that $D$\/ is a squarefree integer.
The change of variable $x \leftrightarrow -x$ shows that
$E_D$ is also isomorphic with $E_{-D}$; this may also be seen from
the Weierstrass equation $y^2 = x^3 - D^2 x$ for~$E_D$.

The arithmetic of the curves $E_D$ has long attracted interest,
both for its connection with the classical ``congruent number''
problem (see \cite[Ch.XVI]{Dickson}; $|D|$ is a ``congruent number''
if and only if $E_D$ has positive rank) and, more recently, as a
paradigmatic example and test case for results and constructions
concerning elliptic curves in general (see for instance \cite{Koblitz}).
The curves $E_D$ have some special properties that make them
more accessible than general elliptic curves over~$\Q$.
They have complex multiplication and are quadratic twists
of the curve $E_1$.  This led to the computation of the sign
of the functional equation of the \hbox{$L$-function} $L(E_D/\Q,s)$:
it depends on $|D|\bmod 8$, and equals $+1$ or $-1$ according as
$|D|$ is in $\{1,2,3\}$ or $\{5,6,7\}$ mod~$8$.  We shall be concerned
with the case of sign~$-1$.

The conjecture of Birch and Swinnerton-Dyer (BSD) predicts that the 
(arithmetic) rank of any elliptic curve~$E$\/ over a number field~$K$,
defined as the \hbox{$\Z$-rank} of its \MW\ group $E(K)$, should equal
the order of vanishing at $s=1$ of $L(E/K,s)$, known as the
``analytic rank'' $\ran(E/K)$.  The BSD conjecture implies the
``BSD parity conjecture'': the arithmetic rank is even or odd
according as the functional equation of $L(E/K,s)$
has sign $+1$ or $-1$.  It would follow that if the sign is $-1$
then $E$\/ always has positive rank.  In our context, where $K=\Q$
and $E=E_D$, this leads to the conjecture that $E_D$ has positive rank
(and thus that $|D|$ is a ``congruent number'')
if $|D|$ is any\footnote{
  We have dropped the hypothesis that $D$\/ be squarefree
  because $c^2 D \equiv D \bmod 8$ for any odd integer~$c$.
  Our integers~$D$\/ are not divisible by~$4$,
  and therefore cannot be of the form $c^2 D$\/ for any even~$c$.
  }
integer of the form $8k+5$, $8k+6$, or $8k+7$.
\subsection{New results and computations}
We prove:
\begin{theorem}
Let $D$\/ be an integer such that
$|D|$ is congruent to $5$, $6$, or~$7$ mod~$8$
and also satisfies $|D|<10^6$.
Then $E_D$ has positive rank over~$\Q$.
\label{thm1}
\end{theorem}
In our ANTS-1 paper~\cite{NDE:A1} we announced such a result
for $|D| < 2\cdot 10^5$.  Our main tool for proving Theorem~\ref{thm1}
is the same: we use the modular parametrization of $E_1$ or $E_2$
to construct a rational point $P_D$ on each $E_D$ from CM points
on modular curves, and usually compute $P_D$ to enough accuracy
to distinguish it from any of the rational torsion points on $E_D$.
Faster computer hardware and new software were both needed to extend
the computation to~$10^6$.  The faster machine made it feasible
to compute $P_D$ for more and larger~$D$.
Cremona's program \textsc{mwrank},
not available when \cite{NDE:A1} was written, found rational points
on the curves $E_D$ on which we could neither distinguish $P_D$ from
a torsion point nor find a rational nontorsion point by direct search.
This happened for $1375$ values of~$|D|$ ---
less than $0.5\%$ of the total, but too many to list here
a rational point on $E_D$ for each such~$D$.  These tables,
and further computational data on the curves $E_D$,
can be found on the Web starting from
\verb_<www.math.harvard.edu/~elkies/compnt.html>_.

Our computations also yield conjectural information on the rank
of~$E_D$: the rank should equal~$1$ if and only if $P_D$ is nontorsion.
In half the cases, those for which $|D|$ or $|D|/2$
is of the form $8k+7$, we obtain this connection from
Kolyvagin's theorem~\cite{Kolyvagin},
which gives the ``if'' direction unconditionally,
and the Gross-Zagier formula~\cite{GZ},
which gives the ``only if'' direction under the BSD conjecture.
Neither Kolyvagin nor Gross-Zagier has been proved to extend
to the remaining cases, when $|D|$ or $|D|/2$ is of the form $8k+5$.
But we expect that similar results do hold in these cases,
and hence that $E_D$ has rank~$1$ if and only if $P_D$ is nontorsion
also when $|D|$ or $|D|/2$ is congruent to~$5$ mod~$8$.
One piece of evidence in this direction is that whenever we found
$P_D$ to be numerically indistinguishable from a torsion point,
the Selmer groups for the \hbox{$2$-isogenies} between $E_D$
and the curve $Dy^2=x^3+4x$ were large enough for $E_D$ to have
arithmetic rank at least~$3$.  We extended the list of curves $E_D$
of conjectural rank $\geq 3$ to $|D|<10^7$ by imposing the
\hbox{$2$-descent} condition from the start
and computing $P_D$ only for those $D$\/ that pass this test.
We find a total of 8740 values of~$|D|$.  The list not only
provides new numerical data on the distribution of quadratic twists
of rank $>1$ with large~$|D|$, but also suggests unexpected biases
in the distribution that favor some congruence classes of~$|D|$'s.

\section{Proof of Theorem~\ref{thm1}}

Let $D$\/ be a squarefree integer such that
$|D|$ is congruent to 5, 6, or~7 mod~8.
Set $K_D=\Q(\sqrt{-|D|}\,)$ if $D$\/ is odd,
and $K_D=\Q(\sqrt{-|D|/2}\,)$ if $D$\/ is even.
Then $K_D$ is an imaginary quadratic field in which
the rational prime~$2$ splits if $D=8k+7$ or $D=16k+14$,
ramifies if $D=8k+5$, and is inert if $D=16k+6$.
A point $P \in E_D(\Q)$ is equivalent
to a \hbox{$K_D$-rational} point $Q$\/ of $E_1$ or~$E_2$
(according as $D$\/ is odd or even)
whose complex conjugate $\overline{Q}$\/ equals $-Q$.
If $Q'$ is any point of $E_1$ or~$E_2$ over~$K_D$
then $Q = Q' - \overline{Q'}$ satisfies $\overline{Q}=-Q$,
and thus amounts to a point of~$E_D$ over~$\Q$.  To prove
Theorem~\ref{thm1} for $E_D$, it will be enough to find
$Q_D \in E_1(K_D)$ or $E_2(K_D)$ and show that the point
$P_D \in E_D(\Q)$ corresponding to $Q_D - \overline{Q_D}$
is not in $(E_D(\Q))_{\rm tors} = E_D[2]$.

We use the modular parametrizations of $E_1$ and $E_2$
by the modular curves $\X_0(32)$ and $\X_0(64)$.
These curves have ``CM~points'' parametrizing cyclic isogenies
of degree~$32$ or~$64$ between elliptic curves
of complex multiplication by the same order in~$K_D$.
If the prime $2$ splits in $K_D$, these points are defined
over the class field of~$K_D$; otherwise they are defined
over a ray class field.  (In the former case, the CM~points
are often called ``Heegner points''; in the latter, \cite{Monsky}
applies the term ``mock Heegner points'', though Birch points out
that Heegner's seminal paper~\cite{Heegner} already used both
kinds of points to construct rational points on~$E_D$,
and the distinction between the two cases was a later development.)
In either case, we obtain a point $Q_D$ defined over~$K_D$
by taking a suitable subset of these CM points,
mapping them to $E_1$ or $E_2$ by the modular parametrization,
and adding their images using the group law of the curve.
See \cite{Birch1,Birch2,Monsky} for more details on these subsets.

Now the key computational point is that the size of each subset
is proportional to the class number of $K_D$,
and thus to $|D|^{1/2}$ when averaged over~$D$.
This is much smaller than the number of terms of the series needed
to numerically estimate $L'(E_D/\Q,1)$, which is on the order of~$D$\/:
as explained for instance in~\cite{BGZ}, for a general elliptic curve
$E/\Q$ of conductor~$N(E)$ it takes $N^{1/2+\epsilon}$ terms
to adequately estimate $L'(E/\Q,1)$, and $N(E_D)=32 D^2$ or $64D^2$
(according as $D$\/ is odd or even) so $N^{1/2}$ is of order~$D$.  
As explained in~\cite{NDE:A1}, the numerical computation
of each CM point as a point on the complex torus $E_1(\C)$ or $E_2(\C)$ 
to within say $10^{-25}$ takes essentially constant time:
find a representative~$\tau$ in a fundamental domain for
the upper half-plane mod $\Gamma_0(32)$ or $\Gamma_0(64)$,
and sum enough terms of a power series for
$\int_\infty^\tau \varphi\, dq/q$
where $\varphi$ is the modular form for $E_1$ or~$E_2$.
Thus it takes time $\Delta^{3/2+\epsilon}$ (and negligible space)
to approximate $Q_D$ for each $|D|<\Delta$.\footnote{
  This computation is particularly efficient in our setting,
  in which $\varphi$ is a CM~form (so most of its coefficients vanish)
  and the normalizers of $\Gamma_0(32)$, $\Gamma_0(64)$ in $\SL_2(\R)$
  can be used to obtain an equivalent $\tau$ with imaginary part
  at least $1/8$ and $\sqrt 3/16$ respectively.  These efficiencies
  represent a considerable practical improvement, though they contribute
  negligible factors $O(\Delta^\epsilon)$ to the asymptotic
  running time of the computation.}

We implemented this computation in {\sc gp}
and ran it for $\Delta=10^6$.
For all but 1375 of the 303979 squarefree values of $|D|<10^6$
congruent to 5, 6, or~7 mod~8, we found that $P_D$ is at distance
at least $10^{-8}$ from the nearest \hbox{$2$-torsion} point of~$E_D$,
and is thus a rational point of infinite order.

For each of the remaining $D$, the point $P_D$ is numerically
indistinguishable (at distance\footnote{
  Here, as in the preceding paragraph, the distance is measured
  on the complex torus representing $E_1(\C)$ or $E_2(\C)$.
  }
at most $10^{-20}$, usually
much less) from a \hbox{$2$-torsion} point.  We believe
that $P_D$ then actually is a torsion point, and thus that
we must find a nontorsion rational point on $E_D$ in some other way.
We did this as follows.  We first searched for rational numbers $x=r/s$
with $|r|,|s|<5\cdot 10^7$ such that $s^4 x = rs(r^2-s^2)$
is $D$\/ times a square for $|D|<10^6$.
This is a reasonable search since we may assume that $\gcd(r,s)=1$,
require that one of the factors $r,s,r+s,r-s$ of $rs(r^2-s^2)$
have squarefree part $f<(4\cdot 10^6)^{1/4}$ and that another
have squarefree part at most $(4\cdot 10^6/f)^{1/3}$,
and loop over those factors.\footnote{
  In fact we removed the factors of~$4$ by using the squarefree parts
  of $(r\pm s)/2$ instead of $r\pm s$ when $r \equiv s \bmod 2$.
  }
This took several hours and found points on all but~$70$
of our $1375$ $E_D$'s.  The remaining curves were handled
by Cremona's \textsc{mwrank} program, which used a {$2$-descent}
on each curve (exploiting its full rational \hbox{$2$-torsion})
to locate a rational point.
This completed the proof of Theorem~\ref{thm1}.

\section{Curves $E_D$ of conjectural rank $\geq 3$}

It might seem surprising that we were able to find a rational point
on each of the 1375 $E_D$'s for which we could not use $P_D$.
Many curves $E_D$, even with $D$\/ well below our upper limit of~$10^6$,
have rank~$1$ but generator much too large to locate with repeated
{$2$-descents} (see for instance \cite{NDE:A1}).  The reason we could
find nontorsion points on the curves $E_D$ with $P_D \in E_D[2]$
is that these are precisely the curves $E_D$ of odd sign
that should have rank at least~$3$, which makes the minimal height
of a non-torsion point much smaller than it can get in the
\hbox{rank-$1$} case.  We explain these connections below,
and then report on our computations that extend to $10^7$
the list of $|D|$ such that $\ran(E_D)$ is odd
and conjecturally at least~$3$.

\subsection{$P_D$ and the rank of $E_D$}

Consider first the cases $D=8k+7$ and $D=16k+14$.
In these cases the prime~$2$, which is the only prime factor
of the conductors of $E_1$ and $E_2$, is split in $K_D$.
Therefore the results of Gross-Zagier~\cite{GZ}
and Kolyvagin~\cite{Kolyvagin} apply to~$P_D$.
The former result gives the canonical height of~$P_D$
as a positive multiple of $L'(E_D,1)$.  Therefore
$\ran(E_D)>1$ if and only if $P_D$ is torsion.
The latter result shows that if $P_D$ is nontorsion
then in fact the arithmetic rank of~$E_D$ also equals~$1$.
Hence any $E_D$ of rank~$3$ or more must be among those for which
we could not distinguish $P_D$ from a torsion point.

The hypotheses of the theorems of Gross-Zagier and Kolyvagin
are not satisfied in the remaining cases $D=8k+5$ and $D=16k+6$.
However, numerical evidence suggests that both theorems generalize 
to these cases as well.  For instance, when $P_D$ is numerically
indistinguishable from a torsion point, $E_D$ seems to have rank~$3$.
For small $|D|$ we readily find three independent points;
for all $|D|$ in the range of our search, $E_D$ and each of the curves
$Dy^2=x^3+4x$ and $Dy^2=x^3-11x\pm14$ isogenous with $E_D$ has
a \hbox{$2$-Selmer} group large enough to accommodate
three independent points.  When $P_D$ is nontorsion
but has small enough height to be recovered from
its real approximation by continued fractions,
we find that it is divisible by~$2$ if and only if the
\hbox{$2$-Selmer} group has rank at least~$5$, indicating
that $E_D$ has either rank $\geq 3$ or nontrivial $\Sha[2]$.
(The former possibility should not occur, and can often be excluded
by \hbox{$2$-descent} on one of the curves isogenous to~$E_D$.)
Both of these observations are consistent with a generalized
Gross-Zagier formula and the conjecture of Birch and Swinnerton-Dyer,
and would be most unlikely to hold if the vanishing of~$P_D$
had no relation with the arithmetic of~$E_D$.
We thus expect that also in these cases $E_D$ should have rank $>1$
if and only if $P_D$ is a torsion point.

\subsection{Rank and minimal nonzero height}

The conjecture of Birch and Swinnerton-Dyer also explains why
curves $E_D$ of rank $\geq 3$ have nontorsion points of height
much smaller than is typical of curves $E_D$ of rank~$1$.
This conjecture relates the regulator of the \MW\ group of $E_D$
with various invariants of the curve, including its real period
and the leading coefficient $L^{(r)}(E_D,1)/r!$ (where $r=\ran(E_D)$).
Now the real period is proportional to $|D|^{-1/2}$.
The leading coefficient is $\ll |D|^{o(1)}$
under the generalized Riemann hypothesis for $L(E_d,s)$,
or even the weaker assumption of the Lindel\"of conjecture
for this family of {$L$-series} (see for instance~\cite[p.713]{IS}).
One expects, and in practice finds, that it is also $\gg |D|^{-o(1)}$
(otherwise $L(E_d,s)$ has zeros $1+it$ for very small positive~$t$).
Thus we expect the regulator to grow as $|D|^{1/2+o(1)}$,
at least if $\Sha$ is small, which should be true for most~$|D|$.
Hence the minimal nonzero height would be at most $|D|^{1/2r}$.
When $r=1$ this grows so fast that already for $|D|<10^4$
there are many curves $E_D$ with generators much too large
to be found by \hbox{$2$-descents}.\footnote{
  The generators can be obtained using the \hbox{CM-point} construction
  in time $|D|^{O(1)}$, but not $|D|^{1/2+o(1)}$
  because $P_D$ must be computed to high accuracy to recognize
  its coordinates as rational numbers from their real approximations.
  Note that in our computations we showed only that $P_D$ is nontorsion
  and did not attempt to determine it explicitly in $E_D(\Q)$.
  }
But for $r\geq 3$ the minimal nonzero height is at most
$|D|^{1/6+o(1)}$, so $|D|$ must grow much larger before
a \hbox{$2$-descent} search becomes infeasible.

\begin{small}
\textit{Remark}\/ on curves curves $E_D$ of even sign:
For such curves we readily determine whether $\ran(E_D)>0$
by using the Waldspurger-Tunnell formula~\cite{Tunnell}
to compute $L(E_D,1)$.  If $L(E_D,1)\neq0$ then $\ran(E_D)=0$
and $E_D$ also has arithmetic rank~$0$ by Kolyvagin
(or even Coates-Wiles~\cite{CW} because $E_D$ has CM).
If $L(E_D,1)=0$ then $\ran(E_D)\geq 2$, and we can prove that
$E_D$ has positive arithmetic rank if we find a nontorsion point.
We expect that the minimal height of such a point is $|D|^{1/4+o(1)}$.
This grows slower than the $|D|^{1/2+o(1)}$ estimate for rank~$1$,
but fast enough that \hbox{$2$-descent} searches fail
for $|D|$ much smaller than our bound of $10^6$.
Even in the odd-rank case that concerns us in this paper,
it is the curves of rank~$3$ that make it hard to extend
Theorem~\ref{thm1} much beyond $\Delta=10^6$: searching for points
on those curves take time roughly $\exp \Delta^{1/6}$,
which eventually swamps the polynomial time $\Delta^{3/2+\epsilon}$
required to find those curves.
\end{small}

\subsection{Computing $E_D$ of conjectural rank $\geq 3$
with $|D|<10^7$}

We extended to $\Delta=10^7$ our search for $P_D$
numerically indistinguishable from torsion points.
These are the curves that we expect to have rank at least~$3$.
Since we do not expect to extend Theorem~\ref{thm1} to $10^7$,
we saved time by requiring that the Selmer groups for the
isogenies between $E_D$ and $Dy^2=x^3+4x$ be large enough
to together accommodate an arithmetic rank of~$3$.
For very large $\Delta$ this is a negligible saving
because most~$D$\/ pass this test.
But it saved a substantial factor in practice for $\Delta=10^7$:
the test eliminated all but $35\%$ of choices of $|D|=16k+14$,
all but $32.1\%$ of $|D|=16k+6$, all but $21.6\%$ of $|D|=8k+5$,
and all but $16.2\%$ of $|D|=8k+7$.
We found a total of 8740 values of~$D$\/
for which $P_D$ appears to be a torsion point.  We expect
that each $P_D$ is in fact torsion and that the corresponding $E_D$
all have rank at least~$3$.  Some $P_D$ might conceivably be a
nontorsion point very close to $E_D[2]$, but this seems quite unlikely;
at any rate no $P_D$ came closer than $10^{-8}$ but far enough
to distinguish from $E_D[2]$.  All the curves probably have rank
exactly~$3$: the smallest $|D|$ known for a curve $E_D$
of rank~$5$ exceeds $4\cdot 10^9$~\cite{Rogers}.  At any rate
none of our curves with $|D|<2\cdot 10^6$ can have rank~$5$:
we applied {\sc mwrank}'s descents-only mode to each of these $E_D$
and the isogenous curves, and in each case obtained
an upper bound of~$3$ or~$4$ on the rank.
Our curves $E_D$ and the isogenous curves include many examples
of conjectural rank~$3$ and nontrivial~$\Sha[2]$.

There are striking disparities in the distribution of our 8740 values
of~$|D|$ among the allowed congruence classes.  The odd classes
$8k+5$ and $8k+7$ account for $2338$ and $2392$ curves $E_D$
of presumed rank~$3$.  But even $|D|$'s are much more plentiful:
there are $4010$ of them, almost as many as in the two odd classes
combined.  This might be explained by the behavior of the
\hbox{$2$-descent}, which depends on the factorization of~$|D|$,
or the fact that we are twisting a different curve:
$E_1$ for odd~$D$\/ and $E_2$ for even~$D$.
But the $4010$ even~$D$\/'s are themselves unequally distributed
between the $16k+6$ and $16k+14$ cases, the former being significantly
more numerous: $2225$ as against $1785$.  (See Figure~1.)
This disparity is much larger than would be predicted by the
\hbox{$2$-descent} test, which in the range $|D|<10^7$ favors
$16k+16$ but only by a factor of $1.09$ whereas $2225$ exceeds $1785$
by almost~25\%.  Note too that the \hbox{$2$-descent} survival rates
would predict a preponderance of $|D|=8k+7$ over $8k+5$, whereas
the two counts are almost identical.  Do these disparities persist
as $\Delta$ increases, and if so why?  Naturally we would also like
to understand the overall distribution of quadratic twists
of rank $\geq 3$, not only for the ``congruent number'' family
but for an arbitrary initial curve in place of~$E_D$.
We hope that the computational data reported here,
and more fully at \verb_<www.math.harvard.edu/~elkies/compnt.html>_,
might suggest reasonable ideas and conjectures in this direction.

\section*{Acknowledgments}

Thanks to Peter Sarnak for the reference~\cite{IS},
and to the referee for several suggestions that improved the paper.
I am grateful to the Packard Foundation for partial financial support
during the preparation of this work.

\vspace*{6ex}

\makefig{Twists with $|D|\equiv 6\bmod 16$ seem to have rank~$3$\\
          much more often than those with $|D|\equiv 14\bmod 16$}
{fig:gaps}{
~\psfig{figure=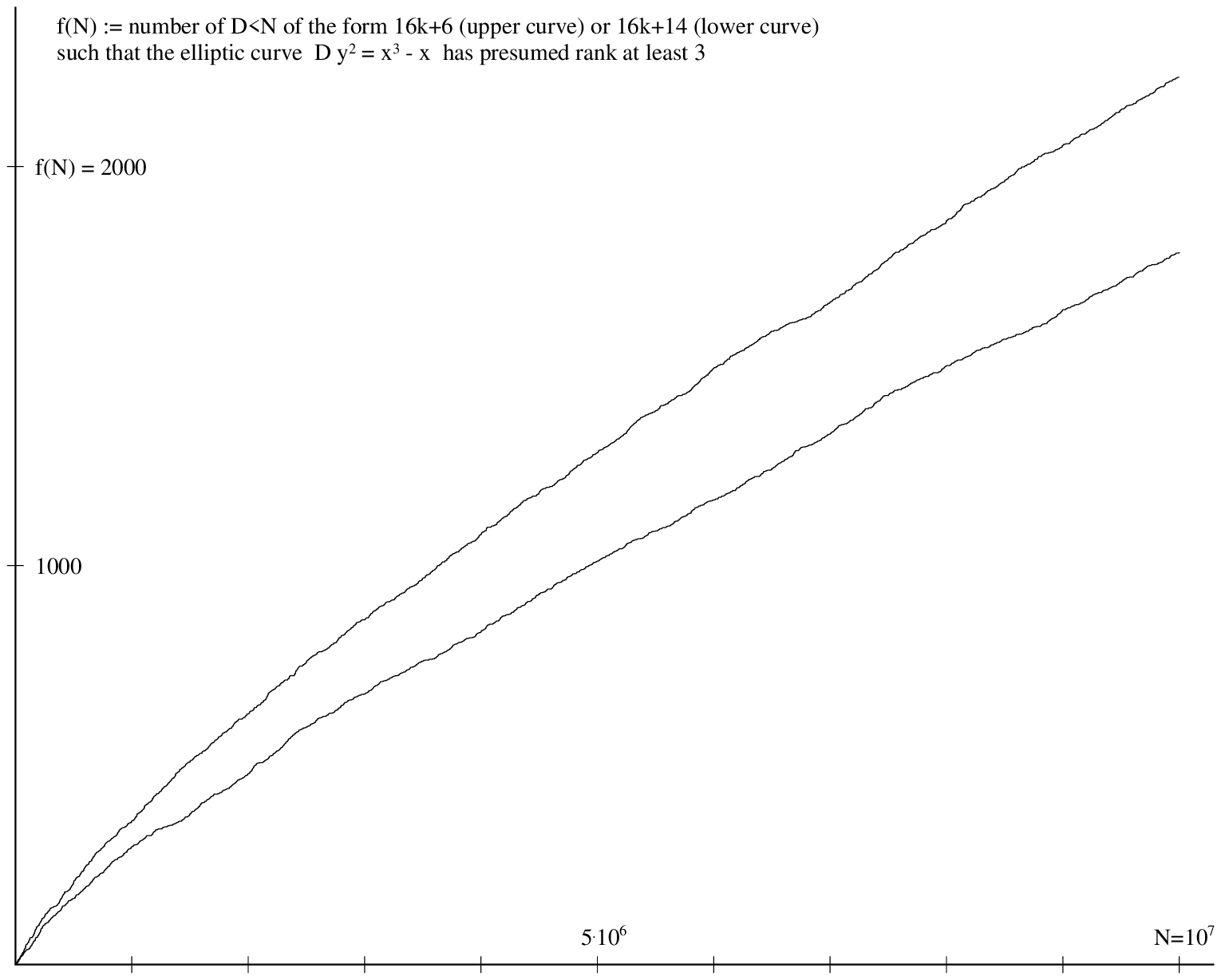,width=\textwidth} \\
}

\end{document}